\documentclass[11pt, twoside, final, reqno]{amsart}
\usepackage{etoolbox,lastpage}
\usepackage{amsmath,amsthm,amscd,amsfonts,amssymb,enumerate}
\usepackage{graphicx}		
\usepackage{color}
\usepackage[colorlinks]{hyperref}
\newtheorem{theorem}{Theorem}[section]

\theoremstyle{definition}

\theoremstyle{remark}

\numberwithin{equation}{section}
\newcommand{\h}{\frac{1}{2}}
 \begin{document}
 \title[{\scriptsize On a new class of series identities}]{On a new class of series identities} 
\author{Arjun K. Rathie }
\address[Arjun K. Rathie]{Department of Mathematics, Vedant College of Engineering \& Technology (Rajasthan Technical    University), Village: Tulsi, Post: Jakhamund, Dist. Bundi, Rajasthan State, India,} 
\email{arjunkumarrathie@gmail.com}
 

 \maketitle

\begin{abstract}The aim of this paper is to provide a new class of series identities in the form of four general results. The results are established with the help of generalizatons of the classical Kummer's summation theorem obtained earlier by Rakha and Rathie. Results obtained earlier by Srivastava, Bailey and Rathie et al. follow special cases of our main findings.
\end{abstract}
 
\section{Introduction and Results Required}
We start with the following two very interesting results involving product of generalized hypergeometric series due to Bailey[1] viz.

\begin{equation}
\begin{aligned}
{}_0F_1& \left[\begin{array}{c} - \\ \rho \end{array}; x\right] \times {}_0F_1\left[\begin{array}{c} - \\ \rho \end{array}; -x\right]\\
&= {}_0F_3\left[\begin{array}{c} - \\ \rho, \frac{1}{2}\rho, \frac{1}{2}\rho +\frac{1}{2}  \end{array}; -\frac{x^2}{4}\right]
\end{aligned}
\end{equation}
and 
\begin{equation}
\begin{aligned}
{}_0F_1& \left[\begin{array}{c} - \\ \rho \end{array}; x\right] \times {}_0F_1\left[\begin{array}{c} - \\ 2-\rho \end{array}; -x\right]\\
&= {}_0F_3\left[\begin{array}{c} - \\ \frac{1}{2}, \frac{1}{2}\rho+\frac{1}{2}, \frac{3}{2}-\frac{1}{2}\rho  \end{array}; -\frac{x^2}{4}\right]\\
& + \frac{2(1-\rho)x}{\rho(2-\rho)} ~{}_0F_3\left[\begin{array}{c} - \\ \frac{3}{2}, \frac{1}{2}\rho+ 1,  2-\frac{1}{2}\rho  \end{array}; -\frac{x^2}{4}\right]
\end{aligned}
\end{equation}
Bailey[1] established these results with the help of the following classical Kummer's summation theorem[2] viz.
\begin{equation}
{}_2F_1\left[\begin{array}{c}a, ~b \\ 1+a-b \end{array}; -1\right] = \frac{\Gamma\left(1+\frac{1}{2}a\right)~\Gamma\left(1+a-b\right)}{\Gamma\left(1+a\right)~\Gamma\left(1+\frac{1}{2}a - b\right)}
\end{equation}
Very recently, Rathie et al.[6] have obtained explicit expressions of\\ 
(i)$\displaystyle  ~{}_0F_1 \left[\begin{array}{c} - \\ \rho \end{array}; x\right] \times {}_0F_1\left[\begin{array}{c} - \\ \rho+i \end{array}; -x\right]$\\
(ii)~ $\displaystyle  {}_0F_1 \left[\begin{array}{c} - \\ \rho \end{array}; x\right] \times {}_0F_1\left[\begin{array}{c} - \\ \rho-i \end{array}; -x\right]$\\
(iii) ~$\displaystyle  {}_0F_1 \left[\begin{array}{c} - \\ \rho \end{array}; x\right] \times {}_0F_1\left[\begin{array}{c} - \\ 2-\rho+i \end{array}; -x\right]$\\
(iv) ~$\displaystyle  {}_0F_1 \left[\begin{array}{c} - \\ \rho \end{array}; x\right] \times {}_0F_1\left[\begin{array}{c} - \\ 2-\rho-i \end{array}; -x\right]$\\
 in the most general form for any $i \in \mathbb{Z}_0$ and provided the natural generalizations of the results (1.1) and (1.2). 
 
 The aim of this paper is to obtain explicit expressions of \\
 (a) ~$\displaystyle  \sum_{m=0}^{\infty} \sum_{n=0}^{\infty}  (-1)^n \frac{\Delta_{m+n}~ x^{m+n}}{(\rho)_m~(\rho+i)_n ~m!~n!}$\\
(b) ~$\displaystyle  \sum_{m=0}^{\infty} \sum_{n=0}^{\infty}  (-1)^n \frac{\Delta_{m+n}~ x^{m+n}}{(\rho)_m~(\rho-i)_n ~m!~n!}$\\
(c) ~$\displaystyle  \sum_{m=0}^{\infty} \sum_{n=0}^{\infty}  (-1)^n \frac{\Delta_{m+n}~ x^{m+n}}{(\rho)_m~(2-\rho+i)_n ~m!~n!}$\\
(d) ~$\displaystyle  \sum_{m=0}^{\infty} \sum_{n=0}^{\infty}  (-1)^n \frac{\Delta_{m+n}~ x^{m+n}}{(\rho)_m~(2-\rho-i)_n ~m!~n!}$\\
in the most general form for any $i \in \mathbb{Z}_0$. Here $\{\Delta_m\}$
is a sequence of arbitrary complex numbers.

The results are derived with the help of the following generalizations of Kummer's summation theorem obtained earlier by Rakha and Rathie[4] for $i \in \mathbb{Z}_0$ viz.
\begin{align}\label{kst1}
	{}_{2}F_{1} \left[ \begin{matrix} a&b\\1+a-b+i\end{matrix};-1\right] 
		&=\frac{2^{-a}\Gamma\left(\frac{1}{2}\right)\Gamma(b-i)\Gamma(1+a-b+i)}{\Gamma(b)\Gamma\left(\frac{1}{2}a-b+\frac{1}{2}i+\frac{1}{2}\right)\Gamma\left(\frac{1}{2}a-b+\frac{1}{2}i+1\right)}\\
		&\times \sum_{r=0}^{i}{i\choose r} (-1)^{r} \frac{\Gamma\left(\frac{1}{2}a-b+\frac{1}{2}i+\frac{1}{2}r+\frac{1}{2}\right)}{ \Gamma\left(\frac{1}{2}a-\frac{1}{2}i+\frac{1}{2}r+\frac{1}{2}\right)}\notag
	\end{align}
and 
	\begin{align}\label{kst2}
	{}_{2}F_{1} \left[ \begin{matrix} a&b\\1+a-b-i\end{matrix};-1\right] 
		&=\frac{2^{-a}\Gamma\left(\frac{1}{2})\Gamma(1+a-b+i\right)}{\Gamma\left(\frac{1}{2}a-b+\frac{1}{2}i+\frac{1}{2}\right)\Gamma\left(\frac{1}{2}a-b+\frac{1}{2}i+1\right)}\\
		&\times \sum_{r=0}^{i}{i\choose r}\frac{\Gamma\left(\frac{1}{2}a-b-\frac{1}{2}i+\frac{1}{2}r+\frac{1}{2}\right)}{ \Gamma\left(\frac{1}{2}a-\frac{1}{2}i+\frac{1}{2}r+\frac{1}{2}\right)}\notag 
	\end{align}
Results obtained earlier by Srivastava[8], Bailey[1] and Rathie et al.[6] follow special cases of our main findings.
\section{Main Results}
The results to be established in this paper are given in the following theorem.
\begin{theorem}
Let $\{\Delta_n\}$ be a bounded sequence of complex numbers. Then for $i \in \mathbb{Z}_0$, the following general results hold true:
\begin{equation}
\begin{aligned}
\sum_{m=0}^{\infty}& \sum_{n=0}^{\infty}  (-1)^n \frac{\Delta_{m+n}~ x^{m+n}}{(\rho)_m~(\rho+i)_n ~m!~n!}\\
&= \sum_{m=0}^{\infty}  \frac{\Delta_m ~x^m}{(\rho)_m ~ m!} ~\frac{2^m~\Gamma\left(\h\right)\Gamma(\rho+i)\Gamma(1-\rho-m-i)}{\Gamma(1-\rho-m)~\Gamma\left( \rho+\h i+\h m - \h\right) \Gamma\left(\rho+\h i + \h m\right)}\\
& \times\left( 
\sum_{r=0}^{i} (-1)^r \binom{i}{r}  \frac{\Gamma\left(\rho+\h m + \h i + \h r-\h\right)}{\Gamma\left( \h r- \h i - \h m + \h \right)}\right)\\
\end{aligned}
\end{equation}
\begin{equation}
\begin{aligned}
\sum_{m=0}^{\infty}& \sum_{n=0}^{\infty}  (-1)^n \frac{\Delta_{m+n}~ x^{m+n}}{(\rho)_m~(\rho-i)_n ~m!~n!}\\
&= \sum_{m=0}^{\infty}  (-1)^n \frac{\Delta_m ~x^m}{(\rho_m ~ m!} ~\frac{2^m~\Gamma\left(\h\right)\Gamma(\rho-i)}{\Gamma\left( \rho-\h i+\h m - \h\right) \Gamma\left(\rho-\h i + \h m\right)}\\
& \times\left( 
\sum_{r=0}^{i}  \binom{i}{r}  \frac{\Gamma\left(\rho+\h m - \h i + \h r-\h\right)}{\Gamma\left( \h r- \h i - \h m + \h \right)}\right)
\end{aligned}
\end{equation}
\begin{equation}
\begin{aligned}
\sum_{m=0}^{\infty}& \sum_{n=0}^{\infty}  (-1)^n \frac{\Delta_{m+n}~ x^{m+n}}{(\rho)_m~(2-\rho+i)_n ~m!~n!}\\
&= \sum_{m=0}^{\infty}  \frac{\Delta_m ~x^m}{(\rho)_m ~ m!} ~\frac{2^{\rho-1+m}~\Gamma\left(\h\right)\Gamma(2-\rho+i)\Gamma(-m-i)}{\Gamma(-m)~\Gamma\left( \h m-\h \rho +\h i+1\right) \Gamma\left(\h m - \h \rho+ \h i \right)}\\
& \times\left( 
\sum_{r=0}^{i} (-1)^r \binom{i}{r}  \frac{\Gamma\left(\h m -\h \rho+ \h i + \h r+1\right)}{\Gamma\left( \h r -\h \rho -  \h m -\h i+ 1 \right)}\right)\\
\end{aligned}
\end{equation}
and    
\begin{equation}
\begin{aligned}
\sum_{m=0}^{\infty}& \sum_{n=0}^{\infty}  (-1)^n \frac{\Delta_{m+n}~ x^{m+n}}{(\rho)_m~(2-\rho-i)_n ~m!~n!}\\
&= \sum_{m=0}^{\infty}  \frac{\Delta_m ~x^m}{(\rho)_m ~ m!} ~\frac{2^{\rho-1+m}~\Gamma\left(\h\right)\Gamma(2-\rho-i)}{\Gamma\left( \h m-\h \rho -\h i+1\right) \Gamma\left(\h m - \h \rho- \h i+\frac{3}{2} \right)}\\
& \times\left( 
\sum_{r=0}^{i}  \binom{i}{r}  \frac{\Gamma\left(\h m -\h \rho- \h i + \h r+1\right)}{\Gamma\left( \h r -\h \rho -  \h m -\h i+ 1 \right)}\right)\\
\end{aligned}
\end{equation}
\end{theorem}	
\subsection*{Derivations :} In order to establish the first result (2.1) asserted in the theorem, we proceed as follows. Denoting the left hand side of (2.1) by $S$, we have
$$S=\sum_{m=0}^{\infty} \sum_{n=0}^{\infty}  (-1)^n \frac{\Delta_{m+n}~ x^{m+n}}{(\rho)_m~(\rho+i)_n ~m!~n!}$$
	Replacing $m$ by $m-n$ and using the result[3, Equ.1, p.56] viz.
	$$\sum_{n=0}^{\infty} \sum_{k=0}^{\infty} A(k, n) = \sum_{n=0}^{\infty} \sum_{k=0}^{n}A(k, n-k)$$
	we have
$$S=\sum_{m=0}^{\infty} \sum_{n=0}^{m}  (-1)^n \frac{\Delta_{m}~ x^{m}}{(\rho)_{m-n}~(\rho+i)_n ~(m-n)!~n!}$$	
Using elementary identities[3, p.58]
$$(\alpha)_{m-n}=\frac{(-1)^n~(\alpha)_m}{(1-\alpha-m)_n}$$
and 
$$(m-n)!=\frac{(-1)^n~m!}{(-m)_n}$$
we have, after some algebra
$$ S=\sum_{m=0}^{\infty}\frac{\Delta_{m}~ x^{m}}{(\rho)_{m}~m!}~\sum_{n=0}^{m}(-1)^n ~\frac{(-m)_n(1-\rho-m)_n}{(\rho+i)_n n!}$$
Summing up the inner series, we have
$$ S= \sum_{m=0}^{\infty} \frac{\Delta_{m}~ x^{m}}{(\rho)_{m}~m!} ~{}_2F_1\left[ \begin{array}{c}-m, 1-\rho-m\\ \rho+i \end{array}; -1\right]$$
We now observe that the series ${}_2F_1$ can be evaluated with the help of the known result (1.4) and we easily arrive at the right hand side of (2.1). This completes the proof of the first result (2.1) asserted in the theorem.

In exactly the same manner, the results (2.2) to (2.4) can be established. So we prefer to omit the details. 
\section{Corollaries}
In this section, we shall mention some of the interesting known as well as new results of our main findings. 
\\
(a) In the result (2.1) or (2.2), if we take $i=0$, we have 
\begin{equation}
\begin{aligned}
\sum_{m=0}^{\infty}& \sum_{n=0}^{\infty}  (-1)^n \frac{\Delta_{m+n}~ x^{m+n}}{(\rho)_m~(\rho)_n ~m!~n!}\\
	&= \sum_{m=0}^{\infty}    \frac{\Delta_{2m}~ (-x^2)^{m}}{(\rho)_m~(\h \rho)_m (\h \rho + \h)_m  ~2^{2m}~m!}
\end{aligned}
\end{equation}
This is a known result due to Srivastava[8]. Further setting $\Delta_m=1 ~(m \in \mathbb{N}_0)$, we at once get the result (1.1) due to Bailey.
\\
(b) In the result (2.3) or (2.4), if we take $i=0$, we have
 \begin{equation}
\begin{aligned}
\sum_{m=0}^{\infty}& \sum_{n=0}^{\infty}  (-1)^n \frac{\Delta_{m+n}~ x^{m+n}}{(\rho)_m~(2-\rho)_n ~m!~n!}\\
	&= \sum_{m=0}^{\infty}    \frac{\Delta_{2m}~ (-x^2)^{m}}{\left(\h\right)_m \left(\h \rho + \h\right)_m~\left(\frac{3}{2}-\h \rho\right)_m  ~2^{2m}~m!}\\
& + \frac{2(1-\rho)x}{\rho(2-\rho)}~\sum_{m=0}^{\infty}  \frac{\Delta_{m}~ (-x^2)^{m}}{\left(\frac{3}{2}\right)_m \left(\h \rho + 1\right)_m~\left(2-\h \rho\right)_m  ~2^{2m}~m!}	
\end{aligned}
\end{equation} 	
which appears to be a new result.

Further setting $\Delta_m=1 ~(m \in \mathbb{N}_0)$, we at once get another  result (1.2) due to Bailey.
\\
(c) In (2.1), if we take $i=0, 1, \cdots, 9$; we get known results recorded in [7].\\
(d) In (2.2), if we take $i=0,1, \cdots, 9$; we get known results recorded in [7].\\
(e)In (2.1) to (2.4), if we get $\Delta_m=1 ~(m \in \mathbb{N}_0)$, we get known results obtained very recently by Rathie et al.[6].

We conclude the paper by remarking that the details about the result presented in this paper together with a large number of special cases(known and new), are given in [5].

\end{document}